\documentclass[12pt,twoside]{article}
\usepackage{amsmath}
\usepackage{amsfonts}
\usepackage{amsthm}
\usepackage{amssymb}
\usepackage{graphicx}
\usepackage{multicol}

\begin{document}
\renewcommand{\baselinestretch}{1.5}
\title{{\normalsize
{\bf  Ribbonness of a stable-ribbon surface-link, I. 
A stably trivial surface-link}}}
\author{{\footnotesize Akio Kawauchi}\\ 
\date{}
{\footnotesize{\it Osaka City University Advanced Mathematical Institute}}\\ 
{\footnotesize{\it Sugimoto, Sumiyoshi-ku, Osaka 558-8585, Japan}}\\ 
{\footnotesize{\it kawauchi@sci.osaka-cu.ac.jp}}}
\maketitle
\vspace{0.25in}
\baselineskip=15pt
\thispagestyle{empty}
\newtheorem{Theorem}{Theorem}[section]
\newtheorem{Conjecture}[Theorem]{Conjecture}
\newtheorem{Lemma}[Theorem]{Lemma}
\newtheorem{Sublemma}[Theorem]{Sublemma}
\newtheorem{Proposition}[Theorem]{Proposition}
\newtheorem{Corollary}[Theorem]{Corollary}
\newtheorem{Claim}[Theorem]{Claim}
\newtheorem{Definition}[Theorem]{Definition}
\newtheorem{Example}[Theorem]{Example}
\begin{abstract}  
There is a question asking whether a handle-irreducible summand of every stable-ribbon surface-link is a unique ribbon surface-link. 
This question for the case of a trivial surface-link  is  affirmatively answered.  
That is,  a handle-irreducible summand of every stably trivial surface-link  
is only a trivial 2-link.  
By combining this result with an old result of  
F. Hosowaka and the author that every surface-knot with infinite cyclic fundamental group is a stably trivial surface-knot, it is concluded  that every surface-knot with infinite cyclic fundamental group is a trivial (i.e., an unknotted) surface-knot. 
\end{abstract}

\phantom{x}



\phantom{x}

\section{Introduction}
A {\it surface-link}  is a closed oriented 
(possibly disconnected) surface $F$  embedded in the 4-space ${\mathbf R}^4$ 
by a smooth (or a piecewise-linear locally flat) embedding. 
When a (possibly disconnected) closed surface $\mathbf F$ is fixed, it is also 
called an $\mathbf F$-{\it link}. If $\mathbf F$ is the disjoint union of 
some copies of the 
2-sphere $S^2$, then it is also called a 2-{\it link}. 
When $\mathbf F$ is connected, it is also called a {\it surface-knot}, 
and a 2-{\it knot} for ${\mathbf F}=S^2$. 

Two surface-links $F$ and $F'$  are {\it equivalent} by an {\it equivalence} $f$  if 
$F$ is sent to $F'$  
orientation-preservingly by an orientation-preserving diffeomorphism (or piecewise-linear homeomorphism) 
$f:{\mathbf R}^4\to {\mathbf R}^4$. The notation $F\cong F'$ is used for equivalent surface-links $F$, $F'$. 
A {\it trivial } surface-link is a surface-link $F$ which is the boundary of the union of mutually disjoint 
handlebodies smoothly embedded in ${\mathbf R}^4$, where a handlebody is a 3-manifold which is a 3-ball, a solid torus or 
a boundary-disk sum of some number of solid tori.  A trivial surface-knot is also called an {\it unknotted } surface-knot. 
A trivial disconnected surface-link is also called an {\it unknotted and unlinked} surface-link. 
For any given closed oriented 
(possibly disconnected) surface $\mathbf F$,  a trivial  $\mathbf F$-link exists uniquely  up to equivalences 
(see \cite{HoK}).  
A {\it ribbon} surface-link   is a surface-link $F$ which is obtained 
from  a trivial $2$-link $O$  
by the  surgery along an embedded 1-handle system (see \cite{K2,K3,K4,K5},
\cite[II]{KSS}). 
 A {\it stabilization} of a  surface-link $F$ is a connected sum 
$F^{\#sT}= F\#  _{k=1}^s  T_k$ of $F$ and a  system of trivial torus-knots 
$T_k\,(k=1,2,\dots,s)$.  
By granting $s=0$, we understand that a surface-link $F$ itself is a stabilization of 
$F$. 
The trivial torus-knot system $T_k\,(k=1,2,...,s)$ is called  the {\it stabilizer} on the stabilization $F^{\#sT}$  of $F$.  

A {\it stable-ribbon} surface-link  is a surface-link $F$ such that a stabilization 
$F^{\#sT}$ of $F$  is a ribbon surface-link.  
For every surface-link $F$, there is a surface-link $F^*$ with minimal total genus 
such that 
$F$ is equivalent to a stabilization of $F^*$.  
The surface-link $F^*$ is called a {\it handle-irreducible summand} of $F$. 
The following question is a central question.

\phantom{x}

\noindent{\bf Question 1.0.}  
A handle-irreducible summand of every stable-ribbon surface-link is a ribbon 
surface-link which is unique up to equivalences ?

\phantom{x}

A {\it stably trivial} surface-link  is a surface-link $F$ such that a stabilization 
of $F$  is a trivial surface-link.
 
In this paper, the following theorem is shown answering  affirmatively this question 
for the case of a stably trivial surface-link.  This question in 
the general case will be answered affirmatively in \cite{K7}.

\phantom{x}

\noindent{\bf Theorem~1.1.}   Any handle-irreducible summand of every stably trivial surface-link is  a trivial  2-link. 

\phantom{x}

The following corollary is directly obtained from Theorem~1.1:

\phantom{x}

\noindent{\bf Corollary~1.2.}  
Every stably trivial surface-link  is a trivial surface-link. 

\phantom{x}

If a surface-knot $F$ has an infinite cyclic fundamental group, then 
$F$ is a TOP-trivial surface-knot, which was shown by Freedman 
for a 2-knot and by \cite{HiK, K11} for a higher genus surface-knot.
In the case of a  piecewise linear surface-knot(equivalent to 
a smooth surface-knot), it is known by \cite{HoK} that 
a stabilization  of the surface-knot $F$ is a trivial  surface-knot, namely  
the surface-knot $F$ is a stably trivial surface-knot. Hence the following corollary is directly obtained from Corollary~1.2 
answering the problem \cite[Problem~1.55(A)]{Kir} on unknotting of a 2-knot positively 
(see \cite{K6} for the surface-link version):

\phantom{x}

\noindent{\bf Corollary~1.3.}  A  surface-knot $F$ is a trivial surface-knot  if 
the fundamental group $\pi_1({\mathbf R}^4\setminus F)$ is an infinite cyclic group.

\phantom{x}

The {\it exterior} of  a surface-knot $F$ is the 4-manifold 
$E=\mbox{cl}({\mathbf R}^4\setminus N(F))$ for a tubular neighborhood $ N(F)$ of $F$ in 
${\mathbf R}^4$. Then the boundary $\partial E$ is a trivial circle bundle over $F$.
 A surface-knot $F$ is {\it of Dehn's type}  if there is a section $F'$ of  $F$ in 
the bundle $\partial E$ such that 
the inclusion $F'\to E$ is homotopic to a constant map.
By  \cite[Corollary~4.2]{HiK}, 
the fundamental group  $\pi_1({\mathbf R}^4\setminus F)$ of a surface-knot $F$ 
of Dehn's type is an infinite cyclic group. 
Thus, we have the following corollary(answering the problem \cite[Problem~1.51)]{Kir} on 
unknotting of a 2-knot of Dehn's type positively): 

\phantom{x}

\noindent{\bf Corollary~1.4.}  A surface-knot of Dehn's type is a trivial surface-knot.

\phantom{x}

Unknotting Conjecture asks whether an $n$-knot $K^n$(i.e., a smooth embedding image 
of the $n$-sphere $S^n$ in the(n+2)-sphere $S^{n+2}$) is unknotted 
(i.e., bounds a smooth $(n+1)$-ball in $S^{n+2}$)  
if and only if the complement $S^{n+2}\backslash K^n$ is homotopy equivalent to $S^1$ 
(see \cite{K1} for example).  
This conjecture was previously 
known to be true for $n>3$ by \cite{Le}, for $n=3$ by \cite{Shan} and for $n=1$ by  \cite{Hom, Papa}. 
The conjecture for $n=2$ was known only in the TOP category  by \cite{Fr}(see also \cite{FQ}).
Corollary~1.3 answers this finally remained smooth unknotting conjecture  affirmatively and hence Unknotting Conjecture can be changed into the following:

\phantom{x}

\noindent{\bf Unknotting Theorem.}  A smooth  $S^n$-knot $K^n$ in $S^{n+2}$ is unknotted 
if and only if the complement $S^{n+2}\backslash K^n$ is homotopy equivalent to $S^1$ for every $n\geq 1$. 

\phantom{x}

A main idea in our argument is to use the surgery of a surface-link on an orthogonal 
2-handle pair, which is much different from the surgery of a surface-link on a single 
2-handle. It is known  that every surface-link $F$ in ${\mathbf R}^4$ is obtained 
from a higher genus trivial surface-knot $F'$ by the surgery of $F'$ on a system of 
mutually disjoint 2-handles, because a handlebody in ${\mathbf R}^4$ is obtained 
from a connected Seifert hypersurface of $F$ by removing mutually disjoint 1-handles  
(see \cite{HoK}). Thus, for example, every 2-twist spun 2-bridge knot in 
\cite{Z} is obtained from a trivial torus-knot $T$ in ${\mathbf R}^4$  
by the surgery of $T$ on a single 2-handle, because it bounds a once-punctured 
lens space as a Seifert hypersurface. 

In Section~2, it is shown that every stably trivial surface-link is a trivial 
surface-link if and only if the uniqueness of an orthogonal 2-handle pair on 
every trivial surface-link holds. 
In Section~3, the uniqueness of every orthogonal 2-handle pair on every surface-link 
is shown, by which Theorem~1.1 is obtained.

\phantom{x}

\section{A triviality condition on a stably trivial surface-link}

A 2-{\it handle} on a surface-link $F$  in ${\mathbf R}^4$  is an embedded 2-handle 
$D\times  I$ on $F$ with $D$ a core disk  such that 
$D\times  I\cap  F =\partial  D\times  I$, where 
$I$ denotes a closed interval containing $0$ and $D\times 0$ is identified with $D$. 
If $D$ is an immersed disk, then call it an {\it immersed 2-handle}.  
Two (possibly immersed) 2-handles $D\times  I$ and $E\times  I$  on $F$ are {\it equivalent} if there is 
an equivalence $f:{\mathbf R}^4\to {\mathbf R}^4$ from $F$ to itself 
such that the restriction  $f|_F:F\to F$ is the identity map 
and  $f(D\times I)=E\times I$.

\begin{figure}[hbtp]
\begin{center}
\includegraphics[width=10cm, height=5cm]{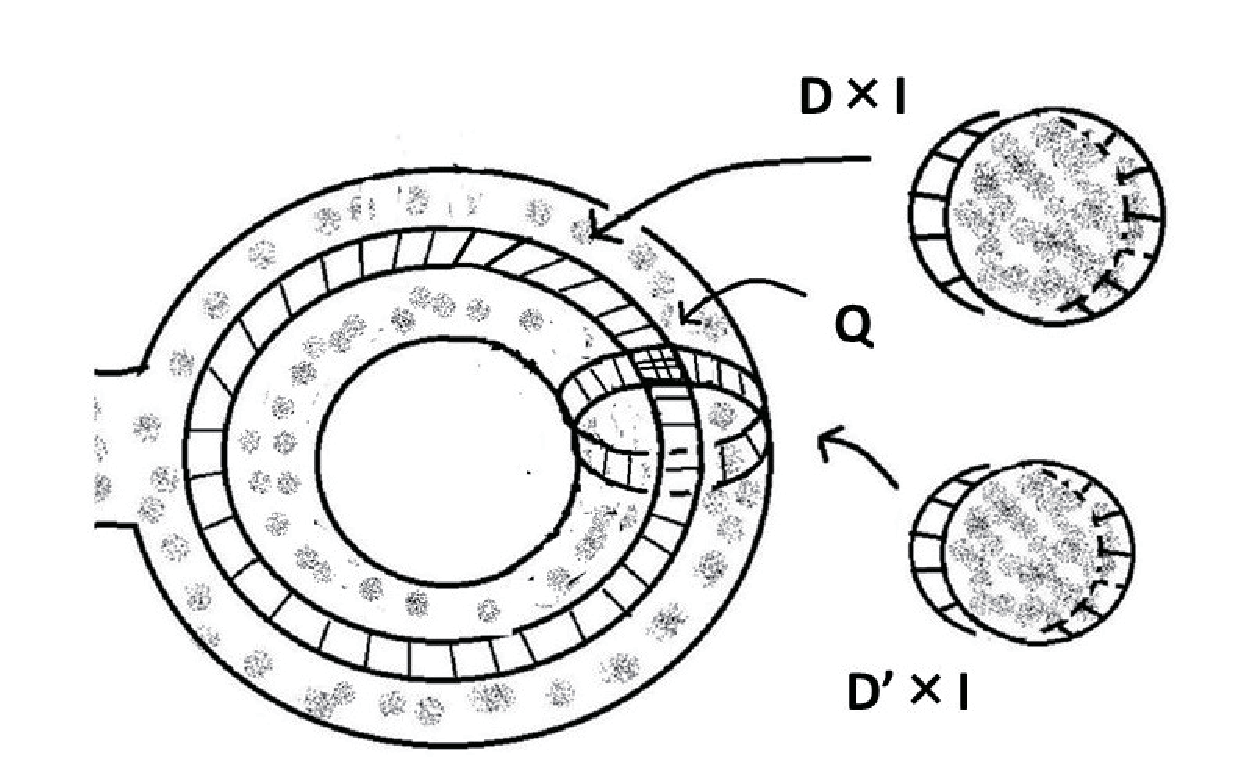}
\end{center}
\caption{An orthogonal 2-handle pair(=: an O2-handle pair) }
\label{fig:O2handles}
\end{figure}

An {\it orthogonal 2-handle pair} (or simply, an {\it O2-handle pair}) on $F$ is a pair   $(D\times I, D' \times  I)$ of 2-handles  
$D\times I$, $D' \times  I$ on  $F$  such that 
\[D\times  I\cap D'  \times  I = 
\partial  D\times  I\cap \partial  D'  \times  I\]
and 
$\partial D\times I$ and  $\partial D' \times  I$ {\it meet orthogonally} on $F$,   
that is, the boundary circles $\partial D$ and  $\partial D'$ meet transversely at one point $p$ and 
the intersection  $\partial D\times I\cap \partial D' \times  I$ 
is homeomorphic to the square $Q=p \times  I\times  I$   
(see  Fig.~\ref{fig:O2handles}).

Let $(D\times I, D'\times I)$  be an O2-handle pair on  a surface-link $F$. 
Let $F(D\times I)$ and  $F(D'\times I)$ be the surface-links obtained from $F$ by the surgeries along $D\times I$ and $D'\times I$, respectively.  
Let $F(D\times I, D'\times I)$ be the surface-link which is  
the union of the plumbed disk 
\[\delta=\delta_{D\times I,D'\times I}
=D\times \partial I\cup Q\cup D'\times partial I\]
and the surface 
\[F^c_{\delta}= \mbox{cl}(F\setminus (\partial D\times I\cup \partial D'\times I).\]  
  
A once-punctured torus $T^o$ in a 3-ball $B$ 
is {\it trivial} if $T^o$ is smoothly and properly embedded  in $B$ which splits 
$B$ into two solid tori.  
A {\it bump} of a surface-link $F$ is a 3-ball $B$ in ${\mathbf R}^4$ with  
$F\cap B=T^o$ a trivial once-punctured torus  in $B$.  
 Let  $F(B)$ be a surface-link $F^c_B\cup \delta_B$   for the surface 
$F^c_B=\mbox{cl}(F\setminus T^o )$ and a disk 
$\delta_B$ in $\partial B$ with $\partial \delta_B  =\partial T^o$, where note that  
$F(B)$ is  uniquely determined up to cellular moves on $\delta_B$ keeping $F^c_B$ fixed. 
Here, a {\it cellular move} of a surface $P$ in ${\mathbf R}^4$ is a surface 
$\tilde P$ in 
${\mathbf R}^4$ such that the complements 
$d=\mbox{cl}(P\setminus P_0)$ and $\tilde d=\mbox{cl}(\tilde P\setminus P_0)$
of the intersection $P_0=P\cap P'$ are disks in the interiors of $P$ and $\tilde P$, 
respectively   
and the union $d\cup \tilde d$ is a 2-sphere bounding a 3-ball smoothly embedded 
in ${\mathbf R}^4$ and not meeting 
$P_0\setminus\partial d=P_0\setminus\partial \tilde d$.  

\begin{figure}[hbtp]
\begin{center}
\includegraphics[width=14cm, height=8cm]{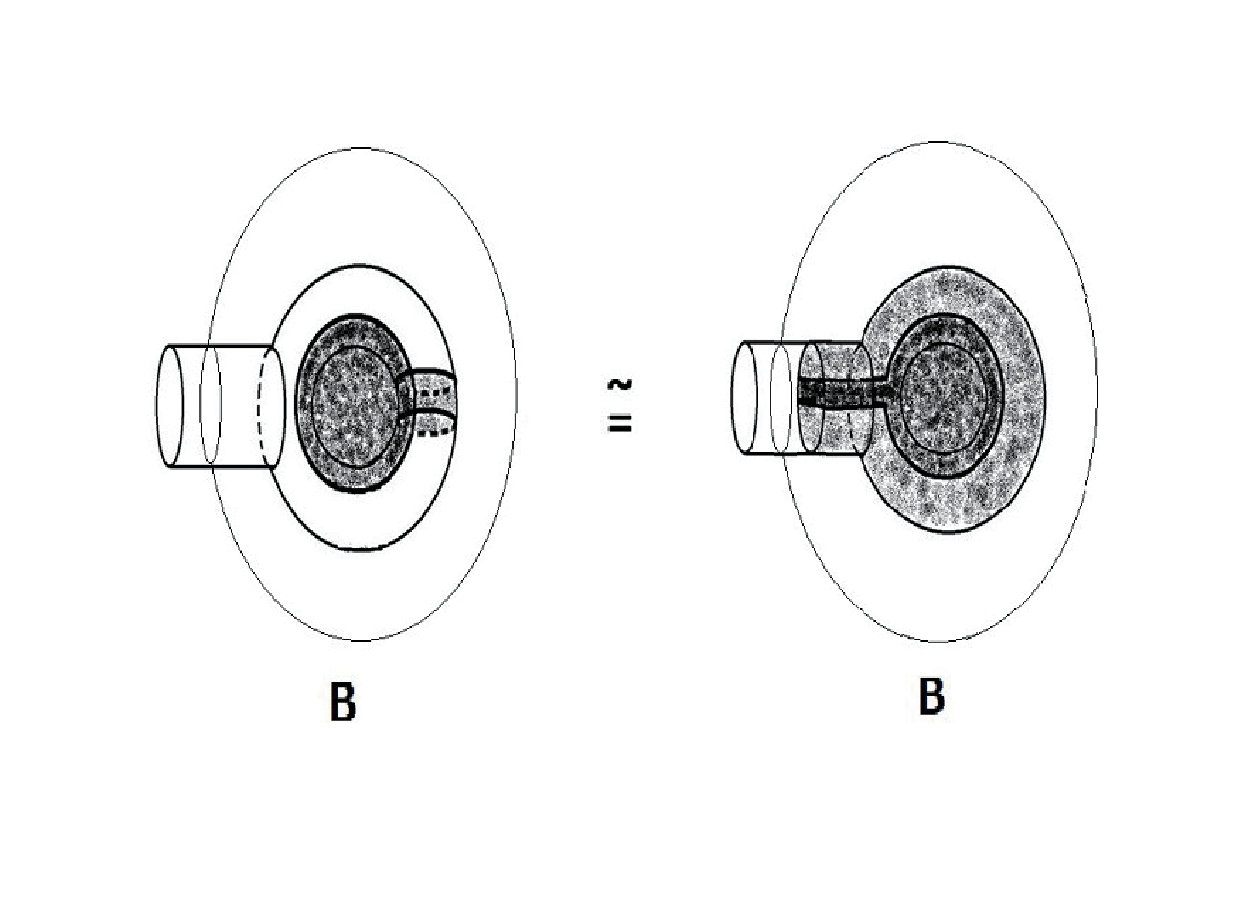}
\end{center}
\caption{ An  associated bump $B$ of a 2-handle union }
\label{fig:bump}
\end{figure}

For an O2-handle pair $(D\times I, D'\times I)$  on  a surface-link $F$, 
let $\Delta=D\times I\cup D'\times I$ is a 3-ball in ${\mathbf R}^4$ called 
the {\it 2-handle union}. Consider  the 3-ball $\Delta$ as a Seifert hypersurface of 
the trivial $S^2$-knot $K=\partial \Delta$ in ${\mathbf R}^4$ to construct a 3-ball 
$B_{\Delta}$ obtained from $\Delta$ by adding an outer boundary collar.  
This 3-ball $B_{\Delta}$ is a bump  of $F$, which we call the 
{\it associated bump}  
of the O2-handle pair  $(D\times I, D'\times I)$. 
When the 3-ball $\Delta$ and a boundary collar of $F^c_{\delta}$ are deformed into 
the 3-space 
${\mathbf R}^3$, this associated bump $B_{\Delta}$ is also considered as a regular neighborhood of $\Delta$ in ${\mathbf R}^3$ (see Fig.~\ref{fig:bump}). 

The following lemma shows that 
giving an O2-handle unordered pair on a surface-link $F$ is the same as giving 
a bump of $F$. 

\phantom{x}

\noindent{\bf Lemma~2.1.} 
An O2-handle unordered pair $(D\times I, D'\times I )$ on a surface-link $F$  
is uniquely constructed from any  given bump $B$ of $F$ in ${\mathbf R}^4$  with 
$F(D\times I, D'\times I)\cong F(B)$. 

\phantom{x}

\noindent{\bf Proof of Lemma~2.1.}
For a bump $B$ of $F$, the set of two solid tori bounded by $T^o =F\cap B$  is 
unique, whose meridian-longitude disk pair is an O2-handle pair. $\square$

\phantom{x}

The following lemma shows the uniqueness of the surgery of  
a surface-link $F$ by an O2-handle pair.

\phantom{x}

\noindent{\bf Lemma~2.2.}
For any O2-handle  pair $(D\times I, D'\times I)$ on any surface-link $F$ and the associated bump $B$, there are equivalences
\[F(B)\cong F(D\times I, D'\times I)  \cong   F(D\times I)   \cong   F(D'\times I).\]
Further, these equivalences are attained by cellular moves  keeping 
$F^c_{\delta}$ fixed. 

\phantom{x}

\noindent{\bf Proof of Lemma~2.2.}  
By definition, we have $F(B)\cong F(D\times I, D'\times I)$. 
The surface-link $F(D\times I, D'\times I)$ is equivalent to $F(D\times I)$ and $F(D'\times I)$ by cellular moves 
on the 3-balls $D'\times I$ and  $D\times I$, respectively. $\square$

\phantom{x}

Two O2-handle pairs $(D\times I, D'\times I)$ and $(E\times I, E'\times I)$ 
on a surface-link $F$ with $\partial D\times I =\partial E\times I$  and 
$\partial D'\times I =\partial E'\times I$ are {\it equivalent} if  
there is an equivalence 
$f: {\mathbf R}^4 \to   {\mathbf R}^4$  from $F$ to itself 
such that the restriction  $f|_F:F\to F$ is the identity map  and   
$f(D\times I)= E\times I$ and  $f(D'\times I)= E'\times I$. 

The following characterization of equivalent O2-handle pairs is useful. 

\phantom{x}

\noindent{\bf Lemma~2.3.} Let $(D\times I, D'\times I)$ and $(E\times I, E'\times I)$ 
be O2-handle pairs on a surface-link $F$ with 
$\partial D\times I =\partial E\times I$  and 
$\partial D'\times I =\partial E'\times I$. 
Let  
\[F(D\times I, D'\times I)=F^c_{\delta}\cup \delta_{D\times I,D'\times I}\quad\mbox{and}\quad  
F(E\times I, E'\times I)=F^c_{\delta}\cup\delta_{E\times I,E'\times I}\]
for the plumbed disks $\delta_{D\times I,D'\times I}$ and 
$\delta_{E\times I,E'\times I}$. 
Then the O2-handle pairs $(D\times I, D'\times I)$ and $(E\times I, E'\times I)$ 
are equivalent if and only if there is an equivalence 
$f: {\mathbf R}^4 \to   {\mathbf R}^4$  from 
$F(D\times I, D'\times I)$ to $F(E\times I, E'\times I)$ 
such that the restriction  $f|_{F^c_{\delta}}:F^c_{\delta}\to F^c_{\delta}$ is the identity map  and   
$f(\delta_{D\times I,D'\times I})= \delta_{E\times I,E'\times I}$. 

\phantom{x}

\noindent{\bf Proof of Lemma~2.3.}  It suffices to show the \lq\lq if\rq\rq part since the \lq\lq only if\rq\rq part is obtained from the definition of equivalent 
O2-handle pairs. Assume that there is an equivalence $f$ from 
$F(D\times I, D'\times I)$ to $F(E\times I, E'\times I)$ 
such that the restriction  $f|_{F^c_{\delta}}:F^c_{\delta}\to F^c_{\delta}$ is the identity map  and   
$f(\delta_{D\times I,D'\times I})= \delta_{E\times I,E'\times I}$. 
The map $f$  is isotopic to a diffeomorphism $f':{\mathbf R}^4\to {\mathbf R}^4$  
sending the associated bump $B_{\Delta(D\times I, D'\times I)}$ of 
$(D\times I, D'\times I)$ to the associated bump 
$B_{\Delta(E\times I, E'\times I)}$ of  $(E\times I, E'\times I)$ by regarding $B_{\Delta(D\times I, D'\times I)}$ 
and $B_{\Delta(E\times I, E'\times I)}$ as  collars of $\delta_{D\times I,D'\times I}$ and $\delta_{E\times I,E'\times I}$, respectively. 
The diffeomorphism $f':{\mathbf R}^4\to {\mathbf R}^4$ is modified into an equivalence 
$f'': {\mathbf R}^4\to {\mathbf R}^4$ from $F$ to itself 
such that the restriction  $f''|_F:F \to F$ is the identity map and
$f''(D\times I)= E\times I$ and  $f''(D'\times I)= E'\times I$. 
Thus, the O2-handle pairs $(D\times I, D'\times I)$ and $(E\times I, E'\times I)$ 
are equivalent.  
$\square$

\phantom{x}

The following corollary is a concrete application of Lemma~2.3.

\phantom{x}

\noindent{\bf Corollary~2.4.} Let $(D\times I, D'\times I)$ and 
$(E\times I, E'\times I)$ be O2-handle pairs on a surface-link $F$ with 
$\partial D\times I =\partial E\times I$  and 
$\partial D'\times I =\partial E'\times I$. 
If the surface-link $F(D\times I, D'\times I)$ 
is obtained from the surface-link $F(E\times I, E'\times I)$ 
by a finite number of cellular moves on $D\times I$, $D'\times I$,  
$E\times I$ and $E'\times I$ keeping $F^c_{\delta}$ fixed, then the O2-handle pairs 
$(D\times I, D'\times I)$ and $(E\times I, E'\times I)$ are equivalent. 

\phantom{x}

\noindent{\bf Proof of Corollary~2.4.} By the assumption, there is an equivalence 
$f: {\mathbf R}^4 \to   {\mathbf R}^4$  from 
$F(D\times I, D'\times I)$ to $F(E\times I, E'\times I)$ 
such that the restriction  $f|_{F^c_{\delta}}:F^c_{\delta}\to F^c_{\delta}$ is 
the identity map  and   
$f(\delta_{D\times I,D'\times I})= \delta_{E\times I,E'\times I}$. 
By Lemma~2.3, the result is obtained.
$\square$

\phantom{x}

A  surface-link $F$ has {\it only unique O2-handle pair} if any two O2-handle pairs 
on $F$ with the same attaching part are equivalent.
A  surface-link not admitting any O2-handle pair  is understood as 
a surface-link with only unique O2-handle pair. 

We have the following characterization on a stably trivial surface-link.

\phantom{x}

\noindent{\bf Lemma~2.5.} The following (1)-(3) are mutually equivalent. 

\medskip 

\noindent{(1)} If a connected sum $F\# T$ of a surface-link $F$ and a trivial 
torus-knot $T$ is  a trivial surface-link, then $F$ is a trivial surface-link.

\medskip 

\noindent{(2)} If $F$ is a trivial surface-link and 
$(D\times I, D'\times I)$ is an O2-handle pair on $F$, then 
$F(D\times I, D'\times I)$ is a trivial surface-link.

\medskip 

\noindent{(3)} Any trivial surface-link has only unique O2-handle pair.

\phantom{x}

\noindent{\bf Proof of Lemma~2.5.} 
$(1) \Rightarrow  (2)$: Let $B$ be the associated bump  of the O2-handle pair 
$(D\times I, D'\times I)$. 
A  4-ball  $A$ obtained by taking a bi-collar  $c(B\times [-1,1])$ of $B$ in 
${\mathbf R}^4$ with $c(B\times 0)=B$ gives a connected sum decomposition 
$F \cong F(D\times I, D'\times I)\# T$.  
By (1),  $F(D\times I, D'\times I)$ is a trivial surface-link.

$(2) \Rightarrow  (3)$:   Let $(D\times I, D'\times I)$ and $(E\times I, E'\times I)$  be  O2-handle pairs with 
$\partial D\times I =\partial E\times I$  and 
$\partial D'\times I =\partial E'\times I$.  
Let 
$F(D\times I, D'\times I)= F^c_{\delta}\cup \delta_{D\times I,D'\times I}$  and  
$F(E\times I, E'\times I)= F^c_{\delta}\cup \delta_{E\times I,E'\times I}$ 
be trivial surface-links for disks $\delta_{D\times I,D'\times I}$ and 
$\delta_{E\times I,E'\times I}$ in the boundaries 
$\partial \Delta(D\times I, D'\times I)$ and $\partial \Delta(E\times I, E'\times I)$ 
of the 2-handle unions 
$\Delta(D\times I, D'\times I)$ and $\Delta(E\times I, E'\times I)$, 
respectively. 
Let $F(D\times I, D'\times I)_0$  
and  $F(E\times I, E'\times I)_0$ be the components of $F(D\times I, D'\times I)$  
and  $F(E\times I, E'\times I)$ containing the loop 
$\partial \delta_{D\times I,D'\times I}=\partial \delta_{E\times I,E'\times I}$, 
respectively, which are made split from the other components in ${\mathbf R}^4$ 
because all the components of every trivial surface-link are split in ${\mathbf R}^4$. 
Since $F(D\times I, D'\times I)_0$ and $F(E\times I, E'\times I)_0$ are trivial 
surface-knots of the same genus, there is 
an equivalence $f: {\mathbf R}^4 \to  {\mathbf R}^4$  sending 
$F(D\times I, D'\times I)_0$ to $F(E\times I, E'\times I)_0$ 
orientation-preservingly and the other components identically.
By a cellular move of $\delta_{D\times I,D'\times I}$ in $F(D\times I, D'\times I)_0$,  this map $f$ is modified to have 
$f(\delta_{D\times I,D'\times I})=\delta_{E\times I,E'\times I}$. 
Further, this map $f$ is modified to send $F^c_{\delta}\cup \delta_{D\times I,D'\times I}$  to 
$F^c_{\delta}\cup \delta_{E\times I,E'\times I}$ 
by sending all the  components except for 
$F(D\times I, D'\times I)_0$ and $F(E\times I, E'\times I)_0$ identically.       
Thus, we have an equivalence $f$ with $f(F^c_{\delta})=F^c_{\delta}$ and 
$f(\delta_{D\times I,D'\times I})= \delta_{E\times I,E'\times I}$. 
By Lemma~2.3, the O2-handle pairs 
$(D\times I, D'\times I)$ and $(E\times I, E'\times I)$ are equivalent.

$(3) \Rightarrow  (1)$: Let $F_i\, (i=0,1,\dots,r)$ be the components of $F$, and 
$F\# T=F_0\#T \cup F_1\cup \dots\cup F_r$ a trivial surface-link. 
Let $V$ be the disjoint union of handlebodies 
$V_i\, (i=0,1,\dots,r)$ in ${\mathbf R}^4$ such that $\partial V_0=F_0\#T$ and 
$\partial V_i=F_i\, (i=1,2,\dots,r)$.

A {\it loop basis} of  $F_0\#T$ of genus $g+1$ is a system of oriented 
simple loop  pairs $(e_j,e'_j)\, (j=0,1,2,\dots,g)$ on $F_0\# T$ representing a basis for $H_1(F_0\# T;{\mathbb Z})$ such that 
$e_j\cap e_{j'}=e'_j\cap e'_{j'}=e_j\cap e'_{j'}=\emptyset$ for all distinct $j,j'$  and 
$e_j\cap e'_j$ is one point  with the intersection 
number $\mbox{Int}(e_j,e'_j)=+1$ in $F_0\#T$ for all  $j$. 
A  loop basis $(e_j,e'_j)\, (j=0,1,2,\dots,g)$ of $F_0\#T$ is {\it spin} 
if the ${\mathbb Z}_2$-quadratic function 
$q:H_1(F_0\# T;{\mathbb Z}_2)\to {\mathbb Z}_2$ associated with the surface-knot 
$F_0\# T$ has  $q(e_j)=q(e'_j)=0$ for all $j$.  
The following result is obtained from \cite[Lemma~2.2]{HiK} where a non-oriented spin loop basis $(e_j,e'_j)\, (j=0,1,2,\dots,g)$ of $F_0\# T$ is constructed.

\phantom{x}

\noindent{(2.5.1)} For a surface-knot $F_0\# T$ of genus $g+1$ in ${\mathbf R}^4$, there is a spin loop basis 
$(e_j,e'_j)\, (j=0,1,2,\dots,g)$ of $F_0\# T$.  In particular, for  a trivial surface-knot $F_0\# T$ bounded by a handlebody $V_0$ in 
${\mathbf R}^4$, every loop basis 
$(e_j,e'_j)\, (j=0,1,2,\dots,g)$ on $\partial V_0$ with  $e'_j\, (j=0, 1,2,\dots,g)$ a meridian loop system of $V_0$
has $q(e'_j)=0$ and either $q(e_j)=0$ or $q(e_j+e'_j)=0$ for all $j$, where 
$e_j+e'_j$ denotes a Dehn twist  of $e_j$ along $e'_j$. 

\phantom{x}

The following result is obtained  from  \cite{Hiro}:

\phantom{x}

\noindent{(2.5.2)}  For any two loop bases 
$(e_j,e'_j)\, (j=0,1,2,\dots,g)$ and $(\tilde e_j,\tilde e'_j)\, (j=0,1,2,\dots,g)$ 
on a trivial genus $g$ surface-knot $F_0\# T$ with 
$q(e_j)=q(\tilde e_j)$ and $q(e'_j)=q(\tilde e'_j)$ for all $j$, there is an 
orientation-preserving diffeomorphism $f:{\mathbf R}^4\to {\mathbf R}^4$ with 
$f(F_0\# T)=F_0\# T$ such that 
$f(e_j)=\tilde e_j$ and $f(e'_j)=\tilde e'_j$ for all $j$.

\phantom{x}

Let $(D\times I, D'\times I)$ be an O2-handle pair on $F\# T$ in ${\mathbf R}^4$ 
attached to $T^o$ such that $(F\# T)(D\times I, D'\times I)\cong F$. 
By (2.5.1), there is a spin loop basis for $F_0\# T$  containing the pair 
$(\partial D, \partial D')$. 
Also, let $(e_i, e'_i)\, (i=0, 1,2,\dots, g)$ be a spin loop basis for 
$F_0\# T$  such that 
$e_0$ bounds a disk $d$ in ${\mathbf R}^4$ 
with $d\cap V= e_0$ and $e'_0$ bounds a meridian disk $d'$ of $V_0$. 
Since the handlebodies 
$V_i\, (i=0,1,\dots,r)$ are splittable in ${\mathbf R}^4$ by \cite{HoK},
we see from (2.5.2) that there is an 
orientation-preserving diffeomorphism $f:{\mathbf R}^4\to {\mathbf R}^4$ with 
$f(F_0\# T)=F_0\# T$ and $f|_{V_i}=1\, (i=1,2,\dots,r)$ such that 
$f(\partial D)=e_0$ and $f(\partial D')=e'_0$. 
A thickening pair $(d\times I, d'\times I)$ of the disk pair $(d,d')$ is an  
O2-handle pair with $(F\# T)(d\times I, d'\times I)$  is a trivial surface-knot. 
Since $(f(D) \times I,f(D')\times I)$ is an O2-handle pair on $F\# T$, we obtain 
from (3) that  
\begin{eqnarray*}
F &\cong& (F\# T)(D \times I,D' \times I)\\
  &\cong& (F\# T)(f(D) \times I,f(D')\times I)\\
  &\cong& (F\# T)(d\times I, d'\times I). 
\end{eqnarray*} 
Thus, $F$ is a trivial surface-link. $\square$

\section{Uniqueness of an orthogonal 2-handle pair}

The following theorem is our main result. 

\phantom{x}

\noindent{\bf Theorem~3.1.} Any (not necessarily trivial) surface-link 
has only unique O2-handle pair.

\phantom{x}

Theorem~1.1 is proved by Theorem~3.1 and Lemma~2.5, which is done as follows:

\phantom{x}

\noindent{\bf Proof of Theorem~1.1.} Let $F$ be a stably trivial link. That is, 
assume that a stabilization $F^{\#sT}= F\#_{k=1}^s T_k$ of $F$ is a trivial link 
for some $s\geq 1$. 
By Theorem~3.1 and Lemma~2.5,  $F\#_{k=1}^{s-1} T_k$ is a trivial surface-link. 
Inductively, $F$  is a surface-link, so that any handle-irreducible summand $F^*$ 
of $F$ is a trivial $S^2$-link. $\square$ 

\phantom{x}

The following lemma is a key lemma to Theorem~3.1.

\phantom{x}

\noindent{\bf Lemma~3.2.} Let  $(D\times I,D'\times I)$ and $(E'\times I,E'\times I)$ 
be O2-handle pairs  on a surface-link $F$ in ${\mathbf R}^4$ with  
$\partial D\times I =\partial E\times I$  and 
$\partial D'\times I =\partial E'\times I$. 
Then there is a 2-handle  $D'_*\times I$ on $F$ 
with $\partial D'_*=\partial D'$ such that the pair
$(E\times I, D'_*\times I)$  is an O2-handle pair on  $F$ 
and the 2-handle $D'_*\times I$ on $F$ is equivalent to the 2-handle $D'\times I$.

\phantom{x}

By assuming Lemma~3.2, the proof of of Theorem~3.1  is  done as follows: 

\phantom{x}

\noindent{\bf Proof of Theorem~3.1.}
Let $(D\times I, D'\times I)$ 
and $(E\times I,E'\times I)$ be  O2-handle pairs on a surface-link  $F$  in 
${\mathbf R}^4$ 
with  $\partial D\times I =\partial E\times I$  and 
$\partial D'\times I =\partial E'\times I$. 
Then there is a 2-handle $D'_*\times I$ on $F$ be a 2-handle on $F$ given by Lemma~3.2 such that 
$(E\times I, D'_*\times I)$  is an O2-handle pair on  $F$ 
and there is an equivalence $f$ from $F$ to itself 
such that the restriction $f|_{F}$ is the identity map on $F$ and  $f(D'_*\times I)= D'\times I$. 
By Lemma~2.2 and Corollary~2.4, the O2-handle pair  $(E\times I, E'\times I)$ on $F$ 
is equivalent to the O2-handle pair $(E\times I, D'_*\times I)$
on $F$, which is equivalent to the O2-handle pair $(f(E)\times I,  D'\times I)$ 
on $F$ and hence  to the  O2-handle pair $(D\times I, D'\times I)$ on $F$. 
Thus,  the O2-handle pair  $(D\times I, D'\times I)$ on $F$  is equivalent 
to an  O2-handle pair $(E\times I, E'\times I)$ on $F$. 
This completes the proof of Theorem~3.1.
$\square$

\phantom{x} 

Throughout the remainder of this section, the proof of Lemma~3.2 is done.

\phantom{x} 

\noindent{\bf Proof of Lemma~3.2.}  
For the core disks $D$, $D'$ $E$ and $E'$ of $D\times I$, $D'\times I$,  
$E\times I$ and $E'\times I$, respectively, assume the following conditions 
(see Fig.~\ref{fig:positions}): 

\begin{figure}[hbtp]
\begin{center}
\includegraphics[width=14cm, height=8.5cm]{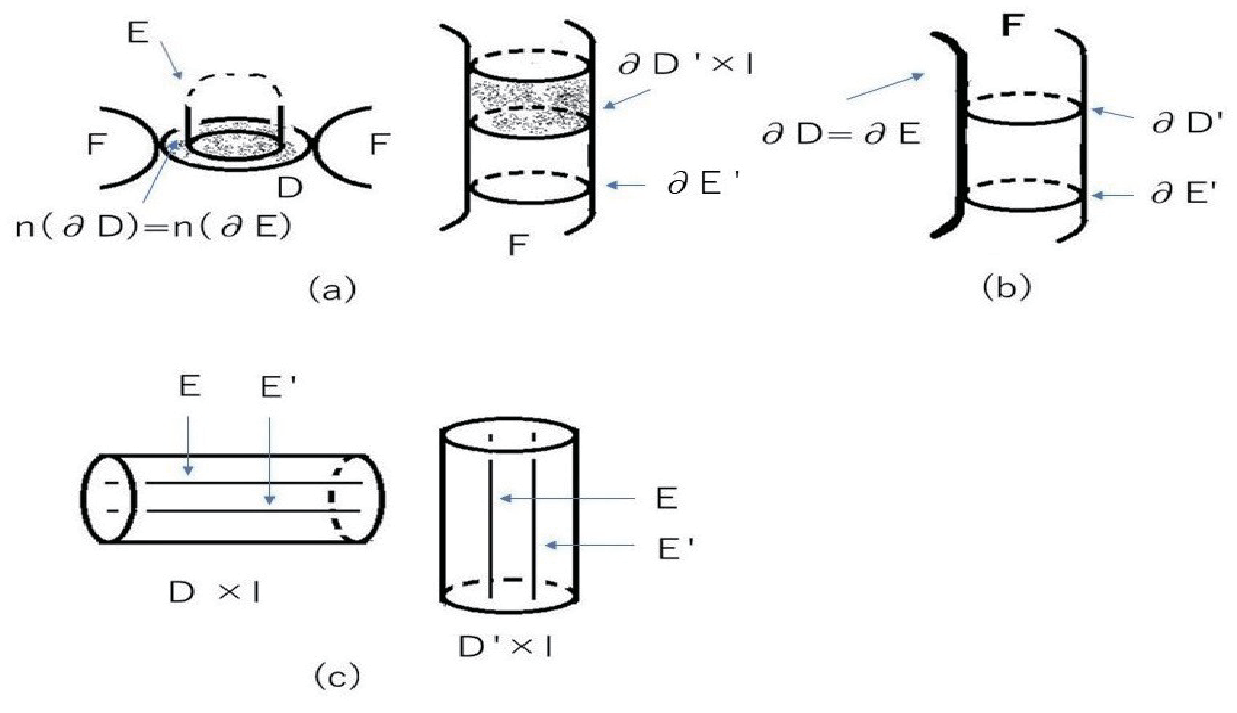}
\end{center}
\caption{Positions among  
the core disks $D$, $D'$, $E$ and $E'$ }
\label{fig:positions}
\end{figure}

\noindent(a) A neighborhood  $n(\partial D)$ of $\partial D$ in $D$ 
coincides with a neighborhood  $n(\partial E)$ of $\partial E$ in $E$  
and  $(\partial D')\times I \cap \partial E'=\emptyset$ by slightly sliding 
$\partial E'$ along $F$,

\medskip

\noindent(b) The disk interiors  $\mbox{Int}D$, $\mbox{Int}D'$, $\mbox{Int}E$ and 
$\mbox{Int}E'$  meet transversely except for the part 
$n(\partial D)=n(\partial E)$ and  
$D\cap D'= \partial D\cap \partial D'=\{p_{D\cap D'}\}$ and 
$E\cap E'=\partial E\cap \partial E'=\{p_{E\cap E'}\}$ for 
distinct points $p_{D\cap D'}$ and  $p_{E\cap E'}$.  

\medskip

\noindent(c) The disk interiors $E\setminus n(\partial E)$ and 
$\mbox{Int}E'$ meet $D\times I$ 
with a finite number of mutually disjoint arcs which are parallel to a fiber $I$ 
of the line bundle $D\times I$ over $D$. 
Similarly, the disk interiors $\mbox{Int}E$ and 
$\mbox{Int}E'$  meet $D'\times I$ 
with a finite number of mutually disjoint arcs which are parallel to a fiber $I$ 
of the line bundle $D'\times I$ over $D'$.

\phantom{x} 

The following operation, called {\it Finger Move Canceling} eliminates an intersection point $x\in \mbox{Int}E\cap \mbox{Int}D'$ by creating a disk $D''$ with  
$\partial D''=\partial D'$ from the disk $D'$.

\phantom{x} 

\noindent{\bf Finger Move Canceling.} 
Let $S$ be a trivial $S^2$-knot  in  ${\mathbf R}^4$ such that 
the 2-sphere $S^2$ is disjoint from $F$ and $D'$ and meets the disk interior 
$\mbox{Int}E$ transversely in just one point $x$. 
Let $y$ be a double point between the disk interiors 
$\mbox{Int}E$ and $\mbox{Int}D'$, 
and $L $ a simple arc in the disk $E$ joining $x$ and $y$ not meeting 
the other double points between $E'$ and  $D$. 
Let $V_{L }$ be a solid tube in  ${\mathbf R}^4$ around the arc $L $  
such that  $V_{L }\cap E=L $ and $V_{L }$ joins  a  disk neighborhood $d_x$ of $x$ in the disk $D'$ and a disk neighborhood $d_y$ of $y$ in the 2-sphere $S$. 
Then a disk $D''$ with $\partial D'' =\partial D'$ and 
$E\cap D''=E\cap D'\setminus \{x\}$ is constructed  
so that  
\[D''=\mbox{cl}(D'\setminus d_x)\cup 
\mbox{cl}(\partial V_{L }\setminus (d_x\cup d_y))\cup 
\mbox{cl}(S \setminus d_y).\]

\phantom{x}

A trivial $S^2$-knot $S$ used in Finger Move Canceling is constructed as follows:

\phantom{x}

\noindent{\bf Claim~3.2.1.} After an isotopic deformation of $F$, $E$ and $E'$ keeping $D$ and $D'$, 
there is a trivial $S^2$-knot $S$ in ${\mathbf R}^4$ such that

\phantom{x}

\noindent(1) $S\cap D=S\cap E=\{ x\}$ for a point $x\in n(\partial D)=n(\partial E)$,

\medskip 

\noindent(2)  $S\cap (F\cup D'\times I \cup E')=\emptyset$, 

\medskip 

\noindent(3) There is a 3-ball $B^S$ in ${\mathbf R}^4$  with $\partial B^S=S$ such that 
$B^S\cap (F\cup D'\times I) =D'$. 

\phantom{x}

By assuming Claim~3.2.1,  let $D'_1$ be a disk parallel to the core disk $D'$ of the 
2-handle $D'\times I$ on the surface-link $F$ such that $D'_1\cap F=\partial D'_1$ 
and  $D'_1\cap (D'\times I)=\emptyset$.  
Let $y$ be  a double point between the disk interiors  
$\mbox{Int}D'_1$ and $\mbox{Int}E$. 
Apply Finger Move Canceling to the trivial $S^2$-knot $S$ in Claim~3.2.1 
along an arc $c$  in $E$ 
from  the point $x$  to  the point $x\in S\cap E$ which avoids the 
double point set $E\cap D'_1\setminus\{y\}$ to obtain a  disk $D'_2$  such that  

\phantom{x} 

\noindent(1) $\partial D'_2 =\partial D'_1$, 

\medskip 

\noindent(2) $E\cap D'_2=(E\cap D'_1)\setminus \{y\}$, and 

\medskip 
 
\noindent(3) $D'_2\cap F=\partial D'_2$ and $D'_2\cap (D'\times I)=\emptyset$.

\phantom{x} 

By continuing this Finger Move Canceling on a trivial $S^2$-knot parallel to $S$, 
a  2-handle  $D'_*\times I$ on $F$ with $\partial D'_*=\partial  D'_1$ such that 
$(E\times I,D'_*\times I)$ is an O2-handle pair on $F$ is obtained. 
The following claim shows that this 2-handle $D'_*\times I$ on the surface-link $F$ is 
a desired 2-handle  in Lemma~3.2. 

\phantom{x} 

\noindent{\bf Claim~3.2.2.} 
The 2-handle $D^*_1\times I$ on $F$ is equivalent to the 2-handle $D'_1\times I$.

\phantom{x} 

This completes the proof of Lemma~3.2 under the assumptions of Claims~3.2.1 and 3.2.2.  

\phantom{x}

The proof of Claim~3.2.1 is done as follows:

\phantom{x}

\noindent{\bf Proof of Claim~3.2.1.}
Let $\Delta$ is the handle union of the O2-handle pair  $(D\times I, D'\times I)$, 
and  $B=B_{\Delta}$  an associated bump of $\Delta$ (see Fig.~\ref{fig:bump}). 
Assume that the bump $B$ is in the 3-space ${\mathbf R}^3$ by an isotopic 
deformation of $B$. 
Let $T^o_B=F\cap B$  be an unknotted once-punctured torus in $B$. 
Let $F^c=\mbox{cl}(F\setminus T^o)$. 
For the sub-surface $T^o_{\Delta}=F\cap \Delta$ of $T^o$, 
the closed complement $A(T^o)=\mbox{cl}(T^o_B\setminus T^o_{\Delta})$ is  an annulus 
bounded by the loops  $o_F=\partial T^o_B=\partial \delta_B=\partial F^c_B$ and 
$o_{\Delta}=\partial T^o_{\Delta}=\partial \delta_{D\times I, D'\times I}$.

Assume that the disk $E$ meets the associated bump $B$ with 
the union of the loop $\partial E$, a set  $J^E_{D\times I}$ of trivial parallel arcs  and a set  $J^E_{D'\times I}$ of trivial parallel arcs such that

\phantom{x}

\noindent(i)  the set  $J^E_{D\times I}$ of  trivial proper parallel arcs in $B$ is obtained by  extending  
the intersection set $\mbox{Int}E\cap (D\times I)$ of trivial parallel arcs in $D\times I$ and

\medskip 

\noindent(ii) the set $J^E_{D'\times I}$ of trivial proper parallel arcs in $B$   is obtained by  extending  
the intersection set  $\mbox{Int}E\cap (D'\times I)$ of trivial parallel arcs in $D'\times I$.

\phantom{x}

Similarly, assume that the disk $E'$ meets the associated bump $B$ with 
the union of the loop $\partial E'$, a set  $J^{E'}_{D\times I}$ of trivial proper  parallel arcs in $B$ and 
a set  $J^{E'}_{D'\times I}$ of trivial proper  parallel arcs in $B$ such that 

\phantom{x}

\noindent(i)${}'$ the set  $J^{E'}_{D\times I}$ of  trivial proper parallel arcs in $B$   is obtained by  extending 
the intersection set $\mbox{Int}E'\cap (D\times I)$ of trivial parallel arcs in 
$D\times I$ and

\medskip 

\noindent(ii)${}'$ the set $J^{E'}_{D'\times I}$ of trivial proper parallel arcs in $B$ is obtained by  extending 
the intersection  $\mbox{Int}E'\cap (D'\times I)$ of trivial parallel arcs in 
$D'\times I$.

\phantom{x}

\begin{figure}[hbtp]
\begin{center}
\includegraphics[width=15cm, height=6cm]{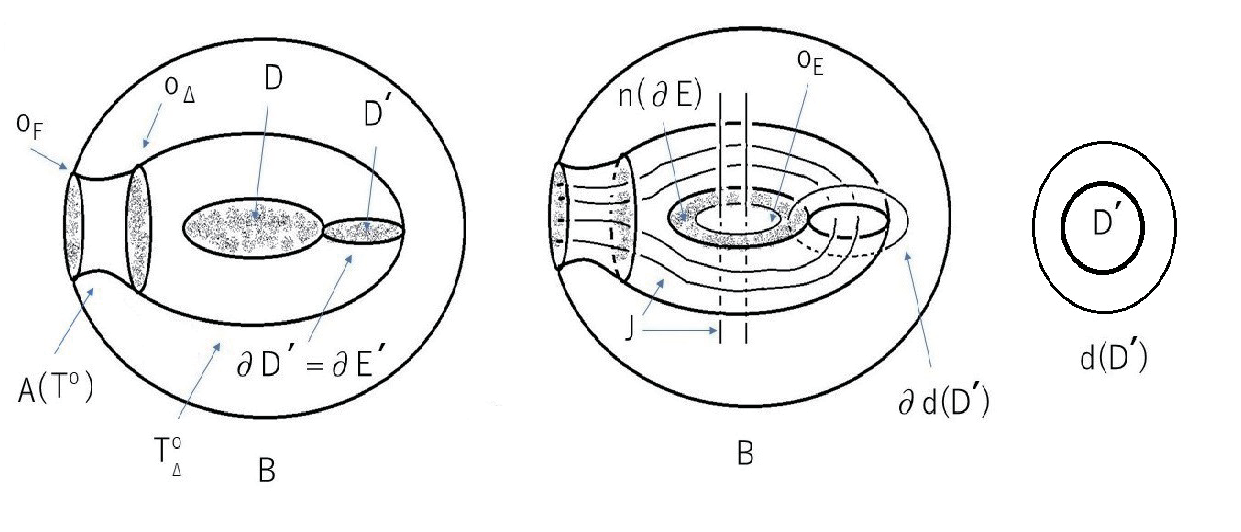}
\end{center}
\caption{A situation of the intersections of the disks $E$ and $E'$ with 
the associated bump $B$ }
\label{fig:Bsections}
\end{figure}

Let 
\[J=J^E_{D\times I}\cup J^E_{D'\times I}\cup J^{E'}_{D\times I}\cup J^{E'}_{D'\times I}.\]
Let $o_E=\partial n(\partial E)\setminus \partial E$. 
Let $d(D')$ be a disk in the associated bump $B$ containing the disk $D'$ 
in the interior such that the link 
$o_E\cup \partial d(D')$ for the boundary loop $\partial d(D')$ is a trivial link in 
$B$ and $\partial d(D')$ transversely meets the disks $E$ and $D$ 
with just one point in the interior of the part $n(\partial D)=n(\partial E)$. 
A situation of the intersections of the disks $E$ and $E'$ with 
the associated bump $B$ of the O2-handle pair  $(D\times I, D'\times I)$ is 
illustrated in Fig.~\ref{fig:Bsections}.

\phantom{x}

\noindent{\bf Notations.}
For a subspace $A$ of ${\mathbf R}^3[0]$ and a subinterval $K$ of ${\mathbf R}$ 
the notation 
\[AK=\{ (x,t)\in {\mathbf R}^4|\, x\in A, t\in K \}\]
is used for a subspace of ${\mathbf R}^4$ as it is used in \cite{KSS}. 
Since the associated bump $B=B_{\Delta}$ of the handle union $\Delta$ of 
the O2-handle pair  $(D\times I, D'\times I)$ is assumed to be  in the 3-space 
${\mathbf R}^3={\mathbf R}^3[0]$, the 4-ball 
\[B[-1,1]\subset {\mathbf R}^3[-1,1]\subset {\mathbf R}^4\]
is a bi-collar of the associated bump of $B$ in the 4-space ${\mathbf R}^4$.
To avoid a confusion, the notation $AK_B$ is used 
for the subspace $AK$ in $B[-1,1]$ defined for a subspace $A$ of $B$ and 
a subinterval $K$ of $[-1,1]$.

\phantom{x}

The following situation may be imposed on the intersection of 
the union $F\cup E \cup E'$ with the 4-ball $B[-1,1]$: 

\phantom{x}

\noindent{\bf (3.2.1.1)}
The surface-knot $F$ and the disks $E$ and $E'$ meet the 4-ball 
$B[-1,1]$ such that 
\[
(F\cup E\cup E')\cap B[t]_B=\left\{
\begin{array}{rl}
(o_{\Delta}\cup J\cup o_E\cup \partial E')[t]_B,& \quad \mbox{for $0< t\leq 1$},\\
(T^o_{\Delta}\cup J \cup n(\partial E))[t]_B,& \quad \mbox{for $t=0$},\\
J[t]_B,& \quad \mbox{for $-1\leq t<0$}. 
\end{array}\right.
\]
  
\phantom{x}

In (3.2.1.1), note that  the annulus $A(T^o)\subset B$ bounded by 
$o_{\Delta}\cup o_F$ is 
deformed into the annulus $o_{\Delta}[0,1]_B\subset B[-1,1]$ 
identifying  $o_{\Delta}\subset B$ with $o_{\Delta}[0]_B\subset B[0]_B$ and 
$o_F\subset B$ with $o_{\Delta}[1]_B\subset B[1]_B$.

\begin{figure}[hbtp]
\begin{center}
\includegraphics[width=14cm, height=8cm]{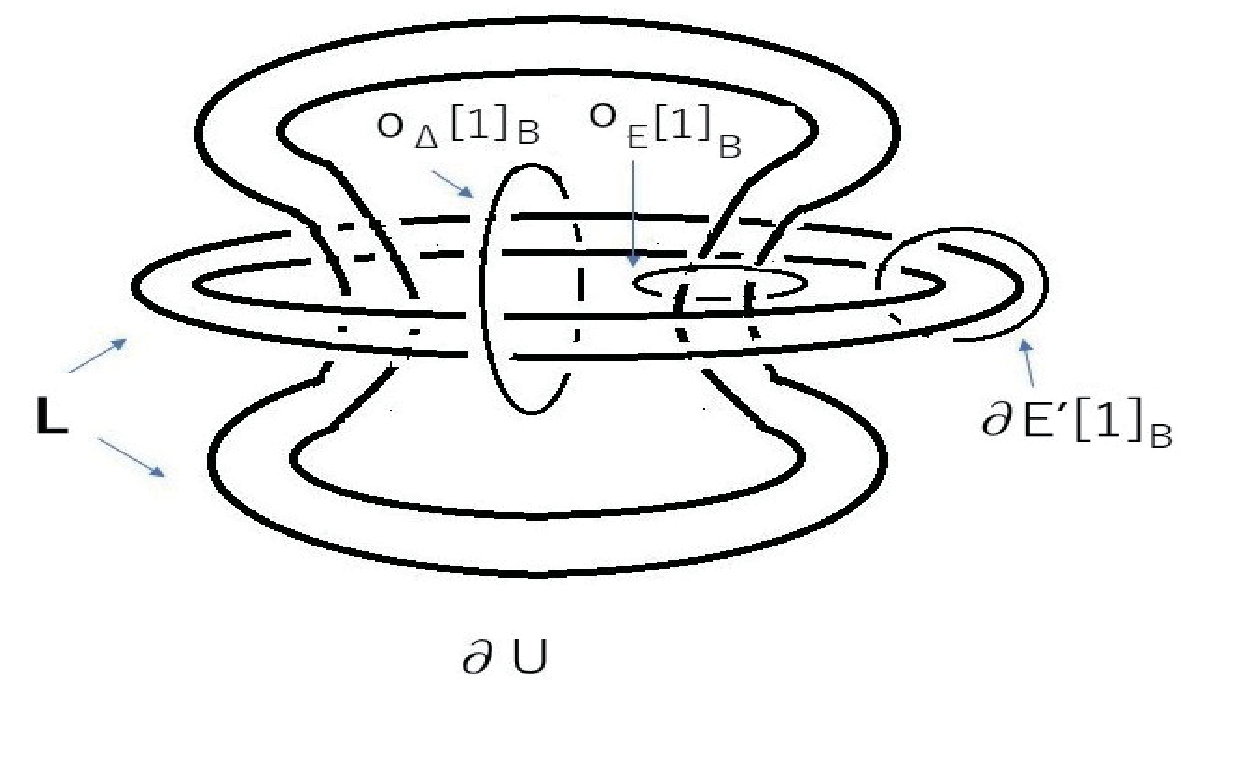}
\end{center}
\caption{The link ${\mathbf L}$ in the 3-sphere $\partial U$}
\label{fig:link}
\end{figure}

Consider the 4-ball $U=\mbox{cl}(\bar{\mathbf R}^4\setminus B[-1,1])$ 
for the one-point-compactification $\bar{\mathbf R}^4$ of the 4-space ${\mathbf R}^4$  
and the proper surfaces
\begin{eqnarray*}
R(F)&=&\mbox{cl}(F\setminus F\cap B[-1,1]),\\
R(E)&=&\mbox{cl}(E\setminus E\cap B[-1,1]),\\
R(E')&=&\mbox{cl}(E'\setminus E'\cap B[-1,1])
\end{eqnarray*}
in the 4-ball $U$. The link 
\[{\mathbf L }=\partial R(F)\cup \partial R(E)\cup \partial R(E')\] 
in the 3-sphere $\partial U=B[-1]_B\cup (\partial B)[-1,1]_B\cup B[1]_B$ 
is illustrated in Fig.~\ref{fig:link}, 
where $\partial R(F)$ and $\partial R(E)\cup \partial R(E')$ are given as follows:
\begin{eqnarray*}
\partial R(F)&=&o_{\Delta}[1]_B\subset \partial U,\\
\partial R(E)\cup \partial R(E')&=& 
o_E[1]_B\cup \partial E'[1]_B\cup {\mathbf L}'\subset \partial U \\
&\phantom{=}&\phantom{XXXXXXX}\mbox{for}\quad 
{\mathbf L}'=J[-1]_B\cup (\partial J)[-1,1]_B\cup J[1]_B.
\end{eqnarray*}

Consider the pair $(U,\partial U)$ as the one-point-compactification 
of the pair of the upper-half 4-space 
\[{\mathbf R}^4_+=\{(x,t)\in{\mathbf R}^3 \times {\mathbf R}
|\, x\in{\mathbf R}^3, t\in {\mathbf R}\}\] 
and the boundary 3-space $\partial {\mathbf R}^4_+={\mathbf R}^3={\mathbf R}^3[0]$.  
The same notations for the proper surface  
$R(F)\cup R(E)\cup R(E')$ in the 4-ball $U$ and the link 
${\mathbf L}=o_{\Delta}[1]_B\cup o_E[1]_B\cup \partial E'[1]_B\cup L'$ in the 
boundary 3-sphere $\partial U$ 
are used for the corresponding proper surface in ${\mathbf R}^4_+$ and the corresponding link in the boundary 3-space ${\mathbf R}^3={\mathbf R}^3[0]$.

By  an argument of \cite{KSS}, a normal form of the surface $R(F)\cup R(E)\cup R(E')$ in ${\mathbf R}^4_+$ is considered 
to obtain the following surface $G$  from  
the surface $R(F)\cup R(E)\cup R(E')$ by an ambient isotopy 
of ${\mathbf R}^4_+$ keeping the boundary ${\mathbf R}^3={\mathbf R}^3[0]$ fixed:

\phantom{x}

\noindent{\bf (3.2.1.2)}  The surface $G$ in ${\mathbf R}^4_+$ is given by 
\[
G\cap {\mathbf R}^3[t]=\left\{
\begin{array}{rl}
\emptyset,&\quad \mbox{for $t>3$}\\
{\mathbf d}({\mathbf O})[t],& \quad \mbox{for $t=3$},\\
{\mathbf O}[t],& \quad \mbox{for $2< t< 3$},\\
({\mathbf L }\cup {\mathbf o}\cup {\mathbf b})[t],& \quad \mbox{for $t=2$},\\
({\mathbf L }\cup {\mathbf o})[t],& \quad \mbox{for $1< t< 2$},\\
({\mathbf L }\cup {\mathbf d})[t],& \quad \mbox{for $t=1$},\\
{\mathbf L }[t],& \quad \mbox{for $0\leq t<1$}, 
\end{array}\right.
\] 
where 
\begin{itemize}
\item ${\mathbf d}$ is a disk system in ${\mathbf R}^3$ disjoint from the link $L $ 
and ${\mathbf o}=\partial {\mathbf d}$, a trivial link, 

\item ${\mathbf b}$ is a band system in ${\mathbf R}^3$ spanning the link 
$L \cup {\mathbf o}$, 

\item ${\mathbf O}$  is a trivial link obtained from the link 
${\mathbf L }\cup {\mathbf o}$  by the surgery along ${\mathbf b}$ 
and  ${\mathbf d}({\mathbf O})$ is a disk system bounding the trivial link 
${\mathbf O}$.
\end{itemize}

\medskip

Let $(d(D')[0], D'[0])$  be the disk pair in ${\mathbf R}^3[0]$ 
corresponding to the disk pair $(d(D')[1]_B, D'[1]_B)$ in the 3-ball 
$B[1]_B\subset\partial U$  obtained from the disk pair $(d(D'), D')$ in $B$. 
Let $(\iota\cdot  d(D')[0],\iota\cdot  D'[0])$ be the disk pair in ${\mathbf R}^3[0]$ 
corresponding to the disk pair $(d(D')[-1]_B, D'[-1]_B)$ in the 3-ball 
$B[-1]_B\subset \partial U$ obtained from the disk pair $(d(D'), D')$ 
in $B$, where  note that the disk pair $(d(D')[-1]_B, D'[-1]_B)$ 
is the image of the disk pair $(d(D')[1]_B, D'[1]_B)$
by  the reflection $\iota $in $B[-1,1]$ sending the point $(x,t)$ to 
the point $(x,-t)$ for $x\in B$ and $t\in[-1,1]$.  

By a replacement to a narrow band  and a band slide on the band system 
${\mathbf b}[2]$ in (3.2.1.2), the following condition cab be imposed: 

\phantom{x}

\noindent{\bf (3.2.1.3)} The band system ${\mathbf b}[2]$ does not meet the disks 
$d(D')[2]$ and $\iota\cdot  d(D')[2]$. 
Thus, for every $t$ with $0\leq t\leq 3$, we have:
\begin{eqnarray*}
d(D')[3]\cap G &=&  d(D')[3]\cap {\mathbf d}({\mathbf O})[3], \\
d(D')[t]\cap G &=& (D'\cap {\mathbf L} )[t], \quad \mbox{for $0\leq t< 3$};\\
\iota\cdot  d(D')[3]\cap G &=& \iota\cdot  d(D')[3]\cap {\mathbf d}({\mathbf O})[3], \\
\iota\cdot  d(D')[t]\cap G &=& (\iota\cdot  D'\cap {\mathbf L} )[t], 
\quad \mbox{for $0\leq t< 3$}.
\end{eqnarray*}

\phantom{x}

Let ${\mathbf p}=D'\cap {\mathbf L}$ be the point system in $B$, and 
${\mathbf p}[0]$  the  point system in  ${\mathbf R}^3[0]$ 
representing the point system ${\mathbf p}[1]_B$ in 3-ball $B[1]_B\subset \partial U$. 
Similarly, let  $\iota\cdot{\mathbf p}[0]$ be the  point system in  
${\mathbf R}^3[0]$   representing the point system ${\mathbf p}[-1]_B$ 
in 3-ball $B[-1]_B\subset \partial U$ which is  $\iota$-reflection  image of 
the point system ${\mathbf p}[1]_B$.

In (3.2.1.3),  the intersection $d(D')[3]\cap {\mathbf d}({\mathbf O})[3]$  
is the disjoint union of an improper arc 
system {\boldmath${\alpha}$}[3] joining the point system ${\mathbf p}[3]$ with a point system 
${\mathbf p}^d[3]$ in the loop $\partial d(D')[3]$  and a proper arc system  
{\boldmath$\beta$}[3] in the disk $d(D')[3]$. 

Similarly, the intersection 
$\iota\cdot  d(D')[3]\cap{\mathbf d}({\mathbf O})[3]$  is 
the disjoint union of an improper arc 
system $\iota\cdot\mbox{\boldmath${\alpha}$}[3]$ joining the point system 
$\iota\cdot{\mathbf p}[3]$  with a point system $\iota\cdot{\mathbf p}^d[3]$ in
the loop $\partial \iota\cdot d(D')[3]$ and a proper arc system  
$\iota\cdot\mbox{\boldmath$\beta$}[3]$ in the disk $\iota\cdot d(D')[3]$. 

Let $\mbox{\boldmath$\beta$}^+[3]$ and 
$\iota \cdot \mbox{\boldmath$\beta$}^+[3]$ 
be slightly extended arc systems of the arc systems {\boldmath$\beta$}[3] and 
$\iota\cdot \mbox{\boldmath$\beta$}[3]$ in 
${\mathbf d}({\mathbf O})[3]$, respectively. 
Let {\boldmath$\gamma$} and $\iota\cdot \mbox{\boldmath$\gamma$}$ be 
the arc systems in 
${\mathbf R}^3[3,4]$ obtained respectively by deforming the extended  arc systems 
$\mbox{\boldmath$\beta$}^+[3]$ and $\iota\cdot \mbox{\boldmath$\beta$}^+[3]$  
as follows:

\phantom{x}

\begin{figure}[hbtp]
\begin{center}
\includegraphics[width=14cm, height=8cm]{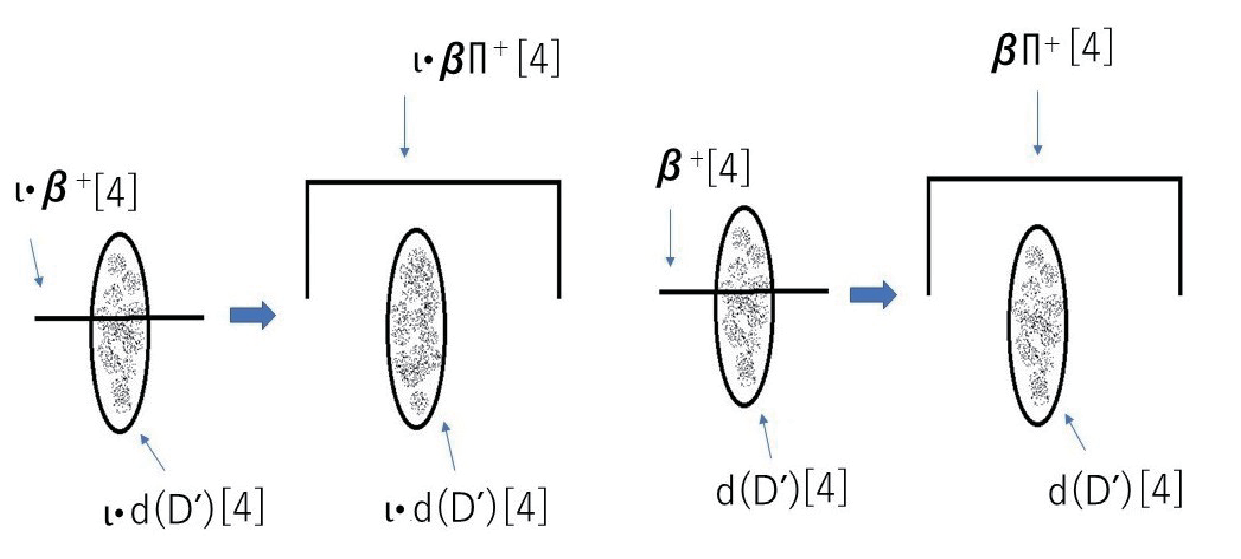}
\end{center}
\caption{The  arc systems $\mbox{\boldmath$\beta$}\sqcap^+[4]$ and 
$\iota\cdot \mbox{\boldmath$\beta$}\sqcap^+[4]$ deformed from 
$\mbox{\boldmath$\beta$}^+[4]$ and 
$\iota\cdot \mbox{\boldmath$\beta$}^+[4]$}
\label{fig:darcsystem}
\end{figure}

\noindent{\bf (3.2.1.4)} For every $t$ with $3\leq t\leq 4$,  
the arc systems $\mbox{\boldmath$\gamma$}$ and $\iota\cdot \mbox{\boldmath$\gamma$}$ 
in ${\mathbf R}^3[3,4]$ are given by 
\[
\mbox{\boldmath$\gamma$}\cap{\mathbf R}^3[t]=\left\{
\begin{array}{rl}
\mbox{\boldmath$\beta$}\sqcap^+[t],& \quad \mbox{for $t=4$},\\
\partial \mbox{\boldmath$\beta$}^+[t],& \quad \mbox{for $3\leq t<4$}, 
\end{array}\right.
\]
where $\mbox{\boldmath$\beta$}\sqcap^+[4]$ 
is an arc system which is deformed from the arc system 
$\mbox{\boldmath$\beta$}^+[4]$ with 
$\partial\mbox{\boldmath$\beta$}\sqcap^+[4]=\partial \mbox{\boldmath$\beta$}^+[4]$ 
and 
$\mbox{\boldmath$\beta$}\sqcap^+[4]\cap d(D')[4]=\emptyset$ 
(see Fig.~\ref{fig:darcsystem}), and 
$$
\iota\cdot \mbox{\boldmath$\gamma$}\cap{\mathbf R}^3[t]
=\left\{
\begin{array}{rl}
\iota\cdot \mbox{\boldmath$\beta$}\sqcap^+[t],& \quad \mbox{for $t=4$},\\
\partial \iota\cdot \mbox{\boldmath$\beta$}^+[t],& 
\quad \mbox{for $3\leq t<4$}, 
\end{array}\right.
$$
where $\iota\cdot \mbox{\boldmath$\beta$}\sqcap^+[4]$ 
is an arc system which is deformed from the arc system 
$\iota\cdot \mbox{\boldmath$\beta$}^+[4]$ 
with 
$\partial \iota\cdot \mbox{\boldmath$\beta$}\sqcap^+[4]=\partial \iota\cdot \mbox{\boldmath$\beta$}^+[4]$ 
and 
$\iota\cdot \mbox{\boldmath$\beta$}\sqcap^+[4]\cap \iota\cdot  d(D')[4]=\emptyset$ 
(see Fig.~\ref{fig:darcsystem}).

\phantom{x}

The deformation from the extended  arc systems $\mbox{\boldmath$\beta$}^+[3]$ and  
$\iota\cdot \mbox{\boldmath$\beta$}^+[3]$ into the arc systems  
{\boldmath$\gamma$}  and  $\iota\cdot \mbox{\boldmath$\gamma$}$ in (3.2.1.4) 
turns the disk system ${\mathbf d}({\mathbf O})[3]$ into a disk system 
${\mathbf d}'({\mathbf O})\subset {\mathbf R}^3[3,4]$ with the intersection 
\[{\mathbf d}''({\mathbf O})[3]={\mathbf d}'({\mathbf O})\cap {\mathbf R}^3[3]\]  
a compact multi-punctured disk system such that  
\[
d(D')[3,4]\cap {\mathbf d}'({\mathbf O})=\mbox{\boldmath${\alpha}$}[3] 
\quad\mbox{and}\quad  
\iota\cdot d(D')[3,4]\cap {\mathbf d}'({\mathbf O})
=\iota\cdot \mbox{\boldmath${\alpha}$}[3].
\]  

Let ${\mathbf q}$ be a point system  in the arc system  
$J^{E'}_{D\times I}\cup J^{E'}_{D'\times I}$ in $B$ which is not in the 
2-handle union $\Delta$. 
Let ${\mathbf a}$ be an arc system in the link ${\mathbf L}$ in $B$ joining the point system  ${\mathbf p}$ with the point system ${\mathbf q}$. 
Let ${\mathbf a}[0]$ and  $\iota\cdot {\mathbf a}[0]$ be the  arc systems in  
${\mathbf R}^3[0]$ representing the arc system ${\mathbf a}[1]_B$ in $B[1]_B$  and 
the arc system ${\mathbf a}[-1]_B$ in $\iota(B[1]_B)=B[-1]_B$, respectively. 
By  a replacement to a narrow band on the band system ${\mathbf b}[2]$ and a band slide, 
assume that the band system ${\mathbf b}[2]$ does not attach to the arc systems 
${\mathbf a}[2]$ and  $\iota\cdot {\mathbf a}[2]$.   
Then the arc systems ${\mathbf a}[3]$ and  $\iota\cdot {\mathbf a}[3]$ are in the boundary of 
the multi-punctured disk system ${\mathbf d}''({\mathbf O})[3]$ 
with  $\partial {\mathbf a}[3]={\mathbf p}[3]\cup {\mathbf q}[3]$ and  
$\partial \iota\cdot {\mathbf a}[3]=\iota\cdot {\mathbf p}[3]\cup 
\iota\cdot  {\mathbf q}[3]$. 

\begin{figure}[hbtp]
\begin{center}
\includegraphics[width=14cm, height=8cm]{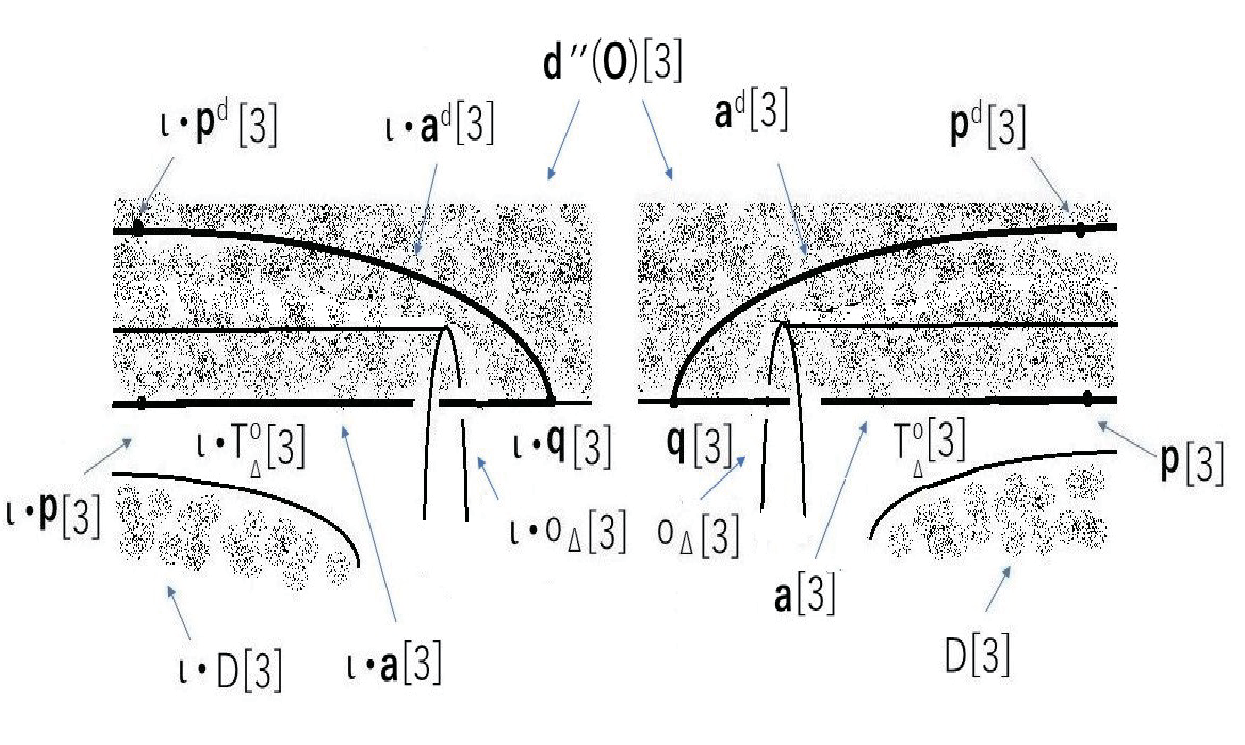}
\end{center}
\caption{Arc systems ${\mathbf a}^d[3]$ and $\iota\cdot {\mathbf a}^d[3]$ }
\label{fig:stararcsystem}
\end{figure}

Let ${\mathbf a}^d[3]$ and $\iota\cdot{\mathbf a}^d[3]$ be arc systems in the 
multi-punctured disk system 
${\mathbf d}''({\mathbf O})[3]$ 
such that $\partial {\mathbf a}^d[3]={\mathbf p}^d[3]\cup {\mathbf q}[3]$ and  
$\partial\iota\cdot{\mathbf a}^d[3] = \iota\cdot{\mathbf p}^d[3]\cup 
\iota\cdot{\mathbf q}[3]$. 
See Fig.~\ref{fig:stararcsystem} for this situation where $T^0_{\Delta}[3]$ and 
$\iota\cdot T^0_{\Delta}[3]$ denote the copies of $T^0_{\Delta}\subset B$ in 
${\mathbf R}^3[3]$ via the copy in $B[1]$ and the reflection image in
$\iota(B[1])=B[-1]$ for the reflection $\iota$ in $B[-1,1]$, respectively.  
 
Let $n({\mathbf a}^d)[3]$ and $n(\iota\cdot{\mathbf a}^d)[3]$ 
be regular neighborhood disk systems of the arc systems 
${\mathbf a}^d[3]$ and $\iota\cdot {\mathbf a}^d[3]$  in 
the multi-punctured disk system ${\mathbf d}''({\mathbf O})[3]$. 

Let ${\mathbf d}^*({\mathbf O})=\mbox{cl}({\mathbf d}'({\mathbf O})
\setminus (n({\mathbf a}^d)[3]
\cup n(\iota\cdot {\mathbf a}^d)[3]))$, 
and ${\mathbf O}^*[t]$ the trivial link obtained from the trivial link 
${\mathbf O}[t]$ by the surgery along 
the disk systems $n({\mathbf a}^d)[t]$ and $n(\iota\cdot{\mathbf a}^d)[t]$ 
for every $t$ with $2<t<3$. 
Also, let ${\mathbf L} ^*[t]$ be the link obtained from the link 
${\mathbf L} [t]$ by surgery along 
the disk systems $n({\mathbf a}^d)[t]$ and $n(\iota\cdot {\mathbf a}^d)[t]$ 
for every $t$ with $1\leq t\leq 2$. 
Then the surface $G^*$ in ${\mathbf R}^4_+$  which is isotopic to $G$ 
by an ambient isotopy keeping  ${\mathbf R}^3[0]$ fixed is given by 
\[
G^*\cap {\mathbf R}^3[t]=\left\{
\begin{array}{rl}
\emptyset,&\quad \mbox{for $t>4$}\\
{\mathbf d}'({\mathbf O})\cap {\mathbf R}^3[t],& \quad \mbox{for $3<t\leq 4$},\\
{\mathbf d}^*({\mathbf O})[t],& \quad \mbox{for $t= 3$},\\
{\mathbf O}^*[t],& \quad \mbox{for $2< t< 3$},\\
({\mathbf L}^*\cup \cap {\mathbf o}\cup {\mathbf b})[t],& \quad \mbox{for $t=2$},\\
({\mathbf L}^*\cup {\mathbf o})[t],& \quad \mbox{for $1< t< 2$},\\
({\mathbf L}^*\cup {\mathbf d})[t],& \quad \mbox{for $t=1$},\\
{\mathbf L}^*[t],& \quad \mbox{for $0\leq t<1$}. 
\end{array}\right.
\]
Let $J^*[1]_B\cup J^*[-1]_B$ be the  arc system in the 3-sphere 
$\partial (B[-1,1])=\partial U$ 
obtained from $J[1]_B\cup J[-1]_B$ by replacing the link ${\mathbf L}[1]_B$ 
with  the link ${\mathbf L}^*[1]_B$ in $\partial (B[-1,1])=\partial U$.   
 
The multi-punctured disk system ${\mathbf d}''({\mathbf O})[3]$ 
is deformed in 
${\mathbf R}^3[3]$ so that $T^o_{\Delta}[0]$ does not meet the 
neighborhood disk systems $n({\mathbf a}^d)[3]$ and 
$n(\iota\cdot {\mathbf a^d})[3]$. 
Then the arc systems $J^*[1]_B$ and  $J^*[-1]_B$ extend to 
the disk system $J^*[-1,1]_B$ in $B[-1,1]$. 

Let $F^*$, $E^*$ and $E^{\prime*}$ be the deformation results of $F$, $E$ and $E'$ 
using $G^*$ and $J^*[-1,1]_B$,  
which are obtained by isotopic deformations on $F$, $E$ and $E'$ keeping $D$ and $D'$ fixed. 
Let $D^S$ and $\iota\cdot  D^S$ be the disks in ${\mathbf R}^3[0,4]$ defined by 
\[
D^S\cap{\mathbf R}^3[t]=\left\{
\begin{array}{rl}
d(D')[t],& \quad \mbox{for $t=4$},\\
\partial d(D')[t],& \quad \mbox{for $0\leq t<4$}, 
\end{array}\right.
\]
\[
\iota\cdot  D^S\cap{\mathbf R}^3[t]=\left\{
\begin{array}{rl}
\iota\cdot  d(D')[t],& \quad \mbox{for $t=4$},\\
\partial \iota\cdot  d(D')[t],& \quad \mbox{for $0\leq t<4$}, 
\end{array}\right.
\]

Let $S$ be the 2-sphere obtained from the disks $D^S$ and $\iota\cdot  D^S$  by connecting  the tube 
$\partial d(D')[-1,1]_B$ in the 4-ball $B[-1,1]$ bounded by 
the loops $\partial D^S$ and $\partial \iota\cdot  D^S$. 
By construction, this 2-sphere $S$ does not meet the surface-link $F^*$ and the disks $D'$, $E^{\prime*}$ 
and meets the disks $D$ and $E^*$ with just one point in the part 
$n(\partial D)=n(\partial E^*)$. 
By construction, there is a 3-ball $B^S$ in ${\mathbf R}^4$  with $\partial B^S=S$ such that 
$B^S\cap (F^*\cup D'\times I) =D'$. 
Thus, $S$ is a desired 2-sphere. 
This completes the proof of Claim~3.2.1.
$\square$ 

\phantom{x} 

The proof of Claim~3.2.2 is done as follows:

\phantom{x}

\noindent{\bf Proof of Claim~3.2.2.} Let $S$ 
be a trivial 2-knot in Claim~3.2.1. 
Let $D'_1\times I$ be a 2-handle  on $F$ with core disk $D'_1$ 
which is disjoint from $D'\times I$. 

Let $D'_2$ be the disk obtained from the disk $D'_1$ and the 2-sphere $S$ 
by taking the surgery along a 1-handle $h$ joining a disk $d'$ in $D'$ 
and a disk $d$ in the $S^2$-knot $S$ and not meeting the interior of the 3-ball $B^3$.   
Let $D'_2\times I$ be the 2-handle on $F$ with $D'_2$ a core disk 
and with $\partial D'_2\times I=\partial D'_1\times I$ which is obtained 
from  the 2-handle $D'_1\times I$ and a collaring $S\times I$ of the trivial 
$S^2$-knot  $S$ and a collaring $h\times I$ of the 1-handle $h$. 
For the bounded surface  $F^c_1=\mbox{cl}(F\setminus \partial D'_1\times I)$,  
the surface-links $F(D'_1\times I)$ and $F(D'_2\times I)$ are given as follows:
\begin{eqnarray*}
F(D'_1\times I)&=& F^c_1\cup D'_1\times \partial I,\\
F(D'_2\times I)&=& F^c_1\cup D'_2\times \partial I.
\end{eqnarray*}

\begin{figure}[hbtp]
\begin{center}
\includegraphics[width=14cm, height=8cm]{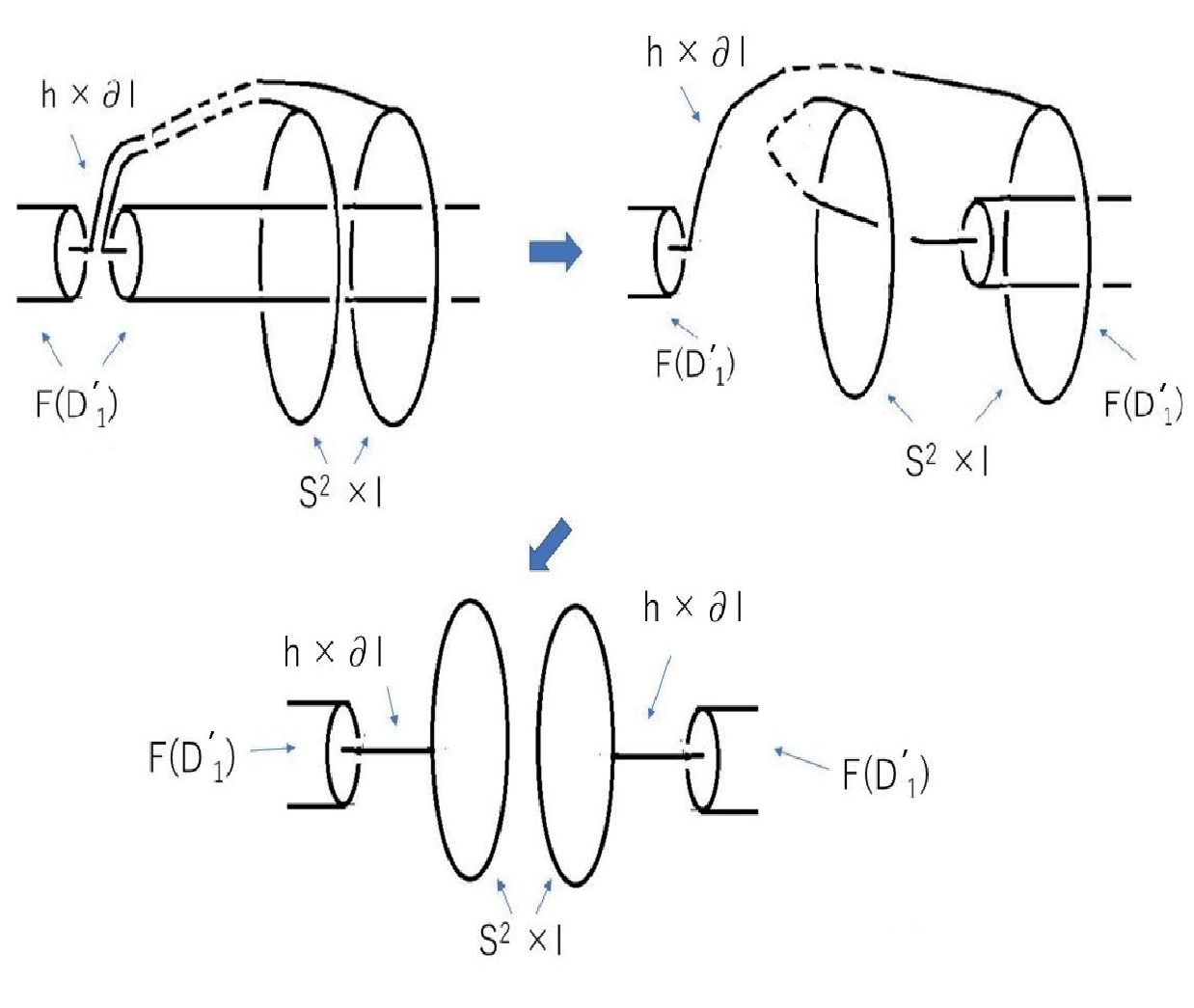}
\end{center}
\caption{An equivalence from the disk $D'_2$ to the disk $D'_1$ }
\label{fig:dashdisk}
\end{figure}

The disk union $D'_2\times \partial I$ is obtained from the disk union 
$D'_1\times \partial I$ by  the surgery along the 1-handle union $h\times \partial I$.  
In Fig~\ref{fig:dashdisk}, it is shown that one 1-handle of the 1-handle union 
$h\times \partial I$ is a self-intersecting 1-handle connecting one disk of 
the disk union $D'_1\times \partial I$ and one 3-ball in 
the 3-ball unions $B^3\times \partial I$ for a collaring $B^3\times I$ of $B^3$. 
This implies that the disk union $D'_2\times \partial I$ 
is deformed into the disk union $D'_1\times \partial I$ by an ambient isotopy of 
${\mathbf R}^4$  keeping the surface $F^c_1$ fixed. 
Thus, there is an equivalence $f:{\mathbf R}^4\to {\mathbf R}^4$ 
from $F(D'_2\times I)$ to $F(D'_1\times I)$ keeping the surface $F^c_1$
identically.

\phantom{x}

The 2-handle $D'_*\times I$ on $F$ constructed by  continuing this operation 
has the property that the pair $(E\times I,D'_*\times I)$ is an O2-handle pair on $F$ 
and there is an equivalence $f:{\mathbf R}^4\to {\mathbf R}^4$ 
from $F(D'_*\times I)$ to $F(D'_1\times I)$ keeping the surface $F^c_1$
identically. 

Let $a'=\partial D\cap D'_1\times I=\partial E\cap D'_*\times I$ be 
the arc parallel to a fiber $I$ 
of the line bundle $\partial D'_1\times I=\partial D'_*\times I$ over  
the circle $\partial D'_1=\partial D'_*$. 
The arc $a'$ attaching to $F(D'_1\times I)$ is $\partial$-relatively isotopic 
to an arc parallel to  $F^c_1$ through the disk $D$. 
Similarly, the arc $a'$ attaching to $F(D'_*\times I)$ is also 
$\partial$-relatively isotopic 
to an arc parallel to  $F^c_1$ through the disk $E$.
This means that the equivalence 
$f$ is isotopically deformed into an equivalence $f'$ 
from $F(D'_*\times I)$ to $F(D'_1\times I)$ keeping the surface $F^c_1$ fixed 
such that $f'(a')=a'$. 
Since the arc $a'$ is regarded as a core of the 1-handle $D'_*\times I$ 
on $F(D'_*\times I)$ and a core of the 1-handle $D'_1\times I$ 
on $F(D'_1\times I)$, the equivalence $f'$ is isotopically deformed into 
an equivalence $f''$ from $F$ to itself such that the restriction 
$f'|_F$ is the identity and  $f''(D'_*\times I)=D'_1\times I$ (see \cite{HoK}).  
This completes the proof of Claim~3.2.2. 
$\square$ 

\phantom{x}

This completes the proof of Lemma~3.2. $\square$ 

\phantom{x}

\noindent{\bf Acknowledgements.} This work was partly supported by 
Osaka City University Advanced
Mathematical Institute (MEXT Joint Usage/Research Center on Mathematics
and Theoretical Physics JPMXP0619217849). 
An idea in this paper together with an idea in \cite{K7} was presented at the meeting 
\lq\lq Differential Topology 19\rq\rq\, held at  Ritsumeikan Tokyo Campus on March 12, 2019, organized by Tetsuya Abe and Motoo Tange. The author would like to thank them
for giving him a talk chance  and the other participants for making some discussions.
The first trial to this research project was done during his stay at Pusan National University in the season of cherry blossoms in full bloom of 2018 spring where 
the author would like to thank Sang Youl Lee and Jieon Kim for kind hospitalities 
during his stay. Since then several improvements of this paper were done and a  
nearly last improvement was done during his stay at Novosibirsk State University for \lq\lq VI Russian-Chinese Conference on Knot Theory and Related Topics\rq\rq held on June 17-21, 2019 where he  would like to thank Nikolay Abrosimov and 
Andrei Vesnin for kind hospitalities during his stay.
The author also would like to thank Seiichi Kamada for taking up a topic 
on a stabilization of a ribbon surface-knot in his lecture \cite{Kamada}, 
which gave him a motivation to consider a stable surface-knot 
and Maggie Miller for giving sharp observations on earlier versions of 
this paper,  Takao Matumoto for telling him where it is difficult to read,  
and  Kengo Kawamura and Masaki Taniguchi talked at OCU Topology Seminars  under  a cooperation of Kouki Sato where some wrong expressions of the paper were pointed out.

\end{document}